\documentclass[11pt]{article}
\usepackage{array, amsmath, amssymb}

\begin{document}
\begin{center}
{\Large\bf Best Simultaneous Approximation of Functions and\\
\vspace{0.3cm} a Generalized Minimax Theorem}\\
\vspace{0.5cm}
{\bf Shinji Tanimoto}\\ 
\vspace{0.5cm}
Department of Mathematics,
University of Kochi,\\
Kochi 780-8515, Japan\footnote{Former affiliation}.  \\
\end{center}
\begin{abstract}
Best simultaneous approximation (BSA) for finitely or infinitely many functions are considered 
under the uniform norm and other important norms. Characterization theorems for a BSA from a finite-dimensional
subspace are obtained by a generalized minimax theorem. From the characterization theorem
a strong unicity theorem is also deduced for a BSA.
\end{abstract}
\vspace{0.7cm}
{\large \bf 1. Introduction} \\  
\\
\indent
Let $\{f_a\}$ be a family of functions obtained in association with each element $a$ in a set $A$.
The purpose is to approximate these functions $\{f_a\}_{a \in A}$  simultaneously from a subspace 
$H$ contained in a function space. In this scetion $X$ is a compact Hausdorff space and
$C(X)$ denotes the set of all real-valued continuous functions on $X$.\\
\indent
In [6] such an approximation problem was considered for real-valued
functions $\{f_a\}_{a \in A}$ defined on $X$. The continuity of
functions themselves is not supposed, but we assume uniform boundedness of the functions. 
For a specified subspace $H$ of finite dimension in $C(X)$,
we say that $f^* \in H$ is a {\it best simultaneous approximation $($BSA$)$} for $\{f_a\}_{a \in A}$ from $H$,
whenever $f^*$ satisfies the inequality
\[ \max_{a \in A, \;x \in X} | f_a(x) - f^*(x) | \le \max_{a \in A, \; x \in X} | f_a(x) - f(x) | 
{\rm ~~for ~all} ~f \in H.\]
In [6] a characterization theorem for a BSA was deduced under the following conditions: 
\begin{itemize} 
\item both functions (of $x$) $\inf_{a \in A}f_a(x)$ and $\sup_{a \in A}f_a(x)$ belong to $C(X)$;
\item for each $x \in X$, the infimum and supremum of $f_a(x)$ are, respectively, attained by some $f_a(x)$.  
\end{itemize}
Moreover, if $H$ is a Haar subspace, a strong unicity theorem for a BSA was obtained from
the characterization theorem (see Section 3). When $X$ is a finite closed interval, an alternation theorem for a BSA 
was also obtained that is similar to the ordinary one (see [1]). \\
\indent
In the next section we consider a BSA problem in a function space $C(X, Y)$ (the set of all continuous functions from
$X$ to $Y$), $Y$ being a normed linear space over the real field $\mathbb{R}$
with norm $\parallel \cdot \parallel$. When a family of functions $\{f_a\}_{a \in A} \subset C(X,Y)$ 
and a finite-dimensional subspace $H \subset C(X,Y)$ are given,
$f^* \in H$ is said to be a BSA to the functions $\{f_a\}_{a \in A}$ from $H$, if the inequality
\[ \max_{a \in A, \;x \in X} \parallel f_a(x) - f^*(x) \parallel \: \le 
\max_{a \in A, \; x \in X} \parallel f_a(x) - f(x) \parallel \]
holds for all $f \in H$. In this setting we will deduce a characterization theorem of a BSA for 
$\{f_a \}_{a \in A}$ that corresponds to the one of [6]. \\
\indent
In Section 3, from this characterization theorem, a strong unicity theorem is derived
in the function space $C(X)$.\\
\indent
In Section 4 we treat another BSA problem in $L^p$-approximation and obtain a
characterization theorem of a BSA for finitely or infinitely many functions $\{f_a\} \subset L^p$. 
These characterization theorems are proved by means of a generalized minimax theorem ([4, Corollary 3.3]). For convenience sake 
we restate it as a lemma. \\
\\
{\bf Lemma 1.}  {\it Let $U$ be an $n$-dimensional, compact convex subset of a Hausdorff topological
vector space, $V$ a compact Hausdorff space, and let $J: U \times V \to \mathbb{R}$ be a jointly
continuous function. An element $u^* \in U$ minimizes $\max_{v \in V}J(u, v)$ over $U$, if and 
only if there exist nonnegative numbers $\lambda_1, \ldots\hspace{-1pt}, \lambda_{n+1}$
with sum one, and $v^*_1, \ldots\hspace{-1pt}, v^*_{n+1} \in V$ such that
\begin{eqnarray}
\sum_{i=1}^{n+1}  \mu_i J(u^*, v_i) \le \sum_{i=1}^{n+1}  \lambda_i J(u^*, v^*_i) 
\le \sum_{i=1}^{n+1}  \lambda_i J(u, v^*_i)  
\end{eqnarray}
holds for all $u \in U$, $v_1, \ldots\hspace{-1pt}, v_{n+1} \in V$, and for all nonnegative numbers $\mu_1, \ldots\hspace{-1pt}, \mu_{n+1}$
with sum one}. \\
\\
\indent
As a useful remark we add that, ignoring all $i$ such that $\lambda_i = 0$ and rearranging the suffix, (1) can be described as
\[ \sum_{i=1}^{n+1}  \mu_i J(u^*, v_i) \le \sum_{i=1}^k  \lambda_i J(u^*, v^*_i) 
\le \sum_{i=1}^k  \lambda_i J(u, v^*_i),  \]
for some $k\, (1 \le k \le n+1)$ with $\sum_{i = 1}^k \lambda_i = 1 ~(\lambda_i > 0)$.\\
\\ 
\\
{\large \bf 2. Characterization  theorem} \\  
\\
\indent
For a compact Hausdorff space $X$ and a normed linear space $Y$ over the real field $\mathbb{R}$
with norm $\parallel \cdot \parallel$, we consider the set $C(X,Y)$ of all continuous functions from
$X$ to $Y$. A family of functions $\{f_a\}_ {a \in A} \subset C(X,Y)$ 
and an $n$-dimensional subspace $H \subset C(X,Y)$ are given, where $n$ is a positive integer.
For $f \in C(X,Y)$ we define the uniform norm of $f$ by
\[  |||f||| = \max_{x \in X} \parallel f(x) \parallel,  \]
and we endow the function space $C(X,Y)$ with this norm. Therefore, a BSA $f^* \in H$ is characterized by
\[ \max_{a \in A} |||f_a - f^*||| \: \le \max_{a \in A} |||f_a - f||| {\rm ~~for~all}~ f \in H. \]
We assume that $A$ is a Hausdorff topological space and impose the two conditions:
\begin{itemize}
\item[\rm(a)] $A$ is compact;  
\item[\rm(b)]  the mapping $A \to C(X, Y)$ defined by $a \mapsto f_a$ is continuous. 
\end{itemize}
\quad \; Now let us introduce the following function
\[ J(f, a, x) = \; \parallel f_a(x) - f(x) \parallel  \]
defined on $H \times A \times X$. It is a jointly continuous function and convex in the argument $f$.
Moreover, $A \times X$ is a compact set with respect to the product topology. Under this setting we
have the following characterization theorem for a BSA. \\
\\
{\bf Theorem 2.}  {\it An element $f^* \in H$ is a BSA to $\{f_a\}$ 
from $H$ if and only if, for some positive integer $k \,(1 \le k \le n + 1)$,
there exist $a_1, \ldots\hspace{-1pt}, a_k \in A$, $x_1, \ldots\hspace{-1pt}, x_k \in X$ and 
positive numbers $\lambda_1, \ldots\hspace{-1pt}, \lambda_k$, whose sum is one, such that}
\begin{itemize}
\item[\rm(i)] 
$\sum_{i=1}^k  \lambda_i \parallel f_{a_i}(x_i) -f^*(x_i)\parallel   \;
\le \sum_{i=1}^k  \lambda_i \parallel f_{a_i}(x_i) -f(x_i)\parallel~~ {\it for~all}~f \in H; $ 
\item[\rm(ii)] $\parallel f_{a_i}(x_i) -f^*(x_i)\parallel \; = |||f_{a_i} -f^*|||
= \max_{a \in A}|||f_a -f^*|||~for ~all~ i \,(1 \le i \le k). $
\end{itemize}
(Proof) Let $f^*$ be a BSA. We define $U = \{f \in H: |||f - f^*||| \le 1 \}$. Then
$U$ is a compact convex set of $H$, since $H$ is finite-dimensional. First we consider
the approximation problem over the set $U$ in place of $H$. Then $f^*$ is also a minimizer
of $\max_{(a, x) \in A \times X}J(f, a, x)$ over $U$. 
Applying Lemma 1 and its remark to this situation, we see that, for some
$k \, (1 \le k \le n+1)$, there exist $(a_1, x_1), \ldots\hspace{-1pt}, (a_k, x_k) \in A \times X$, and
positive numbers $\lambda_1, \ldots\hspace{-1pt}, \lambda_k$ with $\sum_{i=1}^k \lambda_i =1$ such that
the following two inequalities hold: \\
\begin{eqnarray}
\sum_{i=1}^k  \lambda_i J(f^*, a_i, x_i) \le \sum_{i=1}^k  \lambda_i J(f, a_i, x_i) {\rm~ for ~all} ~f \in U;\\
\sum_{i=1}^{n+1}  \mu_i J(f^*, b_i, y_i) \le \sum_{i=1}^k  \lambda_i J(f^*, a_i, x_i) 
\end{eqnarray}
for all $b_1, \ldots\hspace{-1pt}, b_{n+1} \in A$, $y_1, \ldots\hspace{-1pt}, y_{n+1} \in X$ 
and all nonnegative numbers $\mu_1, \ldots\hspace{-1pt}, \mu_{n+1}$ with sum one. \\
\indent
The right-hand side of (2) is a convex function of $f$ and has a local minimum at $f^* \in U$. By a property
of convex functions it follows that it has a global minimum at $f^* \in H$, which implies (i).
Next in (3) putting $\mu_1 = 1$ while other $\mu_i = 0$, and $b_1 = a$ for any $a \in A$, 
we have
$ \parallel f_a(y) -f^*(y) \parallel \; \le \sum_{i=1}^k  \lambda_i \parallel f_{a_i}(x_i) -f^*(x_i)\parallel$
for all $y \in X$, and hence for every $a \in A$ 
\[ |||f_a - f^*||| \le \sum_{i=1}^k  \lambda_i \parallel f_{a_i}(x_i) -f^*(x_i)\parallel. \]  
This shows that
$\max_{a \in A}|||f_a - f^*||| \le \sum_{i=1}^k  \lambda_i \parallel f_{a_i}(x_i) -f^*(x_i)\parallel$.  
Using $\sum_{i = 1}^k \lambda_i = 1$ $(\lambda_i > 0)$, we conclude that
\[\max_{a \in A}|||f_a - f^*||| \le \sum_{i=1}^k  \lambda_i \parallel f_{a_i}(x_i) -f^*(x_i)\parallel \,
\le \max_{a \in A}|||f_a - f^*|||, \]
which implies (ii). \\
\indent
Conversely, suppose that $f^* \in H$ satisfies conditions (i) and (ii) for 
$a_i$'s in $A$, $x_i$'s of $X$ and positive numbers 
$\lambda_i$'s such that $\sum_{i = 1}^k \lambda_i = 1$. Then these conditions imply that, for any $f \in H$, 
\begin{eqnarray*}
\max_{a \in A}|||f_a - f^*||| = \sum_{i=1}^k  \lambda_i \parallel f_{a_i}(x_i) -f^*(x_i)\parallel \, 
\le \sum_{i=1}^k  \lambda_i \parallel f_{a_i}(x_i) -f(x_i)\parallel \, \le \max_{a \in A}|||f_a - f|||,
\end{eqnarray*}
showing that $f^*$ becomes a BSA. This completes the proof. \\
\\
\indent
Next we consider the case where $A$ is a finite set, as discussed in [5].  
Let $g_1, \ldots\hspace{-1pt}, g_{\ell} \in C(X, Y)$ be given. In order to consider BSA to $\{g_j\}$, 
we introduce a compact set 
\[ A = \Big\{ a = (\alpha_1, \ldots\hspace{-1pt}, \alpha_{\ell}) : \sum_{j = 1}^{\ell}\alpha_j = 1, ~\alpha_j \ge 0~(1 \le j \le \ell) \Big\}.\]
For each $a \in A$ we set $g_a = \sum^{\ell}_{j = 1} \alpha_j g_j$.
For $f \in H$, an $n$-dimensional subspace of $C(X, Y)$, we have, using the convexity of norm
\[ \max_{1 \le j \le \ell} ||| g_j - f ||| \le \max_{a \in A}||| g_a - f ||| 
= \max_{a \in A}||| \sum^{\ell}_{1 \le j \le \ell} \alpha_j (g_j -f) ||| 
\le \max_{1 \le j \le \ell} ||| g_j - f |||. \]
Thus our approximation problem is reduced to simultaneously approximate $\{g_a\}~(a \in A)$
from $H$. Then as a special case of Theorem 2 follows the characterization theorem in [5]. \\
\\
{\large \bf 3. Strong unicity theorem} \\  
\\
\indent
Suppose that the norm of $Y$ is defined by means of an inner product $\langle~, ~\rangle$ so that
${\parallel y \parallel}^2 = \langle y, y\rangle$ for $y \in Y$. The condition (i) of Theorem 2
is equivalent to the assertion; the function of a real variable $t$ 
\[ \sum_{i=1}^k  \lambda_i \parallel f_{a_i}(x_i) -f^*(x_i) + tf(x_i) \parallel \]
attains the minimum at $t = 0$ for every $f \in H$. By a simple calculation we obtain the next corollary,
using the inner product. \\
\\
{\bf Corollary 3.}  {\it Let $Y$ be an inner product space. An element $f^* \in H$ is a BSA to $\{f_a\}$ 
from $H$ if and only if, for some positive integer $k \,(1 \le k \le n + 1)$,
there exist $a_1, \ldots\hspace{-1pt}, a_k \in A$, $x_1, \ldots\hspace{-1pt}, x_k \in X$ and 
positive numbers $\lambda_1, \ldots\hspace{-1pt}, \lambda_k$, whose sum is one, such that}\\
\lefteqn{~~{\rm(i')}~~~\sum_{i=1}^k  \lambda_i \langle f_{a_i}(x_i) -f^*(x_i), f(x_i) \rangle = 0~~for~all~f \in H;}\\
\lefteqn{~~{\rm(ii)} ~~\parallel f_{a_i}(x_i) -f^*(x_i)\parallel \; = |||f_{a_i} -f^*|||
= \max_{a \in A}|||f_a -f^*|||~for ~all~ i \,(1 \le i \le k). }
\\
\\
\indent
In what follows we consider the case of $Y = \mathbb{R}$ and so we deal with the function space $C(X)$ as in Section 1.
Referring to [2, p.91] or [4, Section 5], an $n$-dimensional subspace $H \subset C(X)$ is called a 
Haar subspace if, for $n$ distinct elements $x_1, \ldots\hspace{-1pt}, x_n \in X$ and for
$n$ arbitrary numbers $r_1, \ldots\hspace{-1pt}, r_n \in \mathbb{R}$, there exists a unique 
$f \in H$ such that $f(x_k) = r_k \;(1 \le k \le n)$. \\
\indent
A one-dimensional Haar subspace is obviously
spanned by any function that does not vanish in $X$. A two-dimensional Haar subspace is spanned
by every pair of functions $f, \;g \in C(X)$ satisfying $f(x)g(y) \ne f(y)g(x)$ whenever $x \ne y$,
and so on (by linear algebra). 
Here the uniform norm of $C(X)$ is defined by
\[ \parallel f \parallel = \max_{x \in X} |f(x)| \]
for $f \in C(X)$. Now we prove the following strong unicity theorem.\\
\\
{\bf Theorem 4.}  {\it Suppose that the compact set $X$ contains at least $n+1$ elements, where
$n \ge 1$. Let a BSA $f^* \in H$ to $\{f_a\}_{a \in A}$ from an $n$-dimensional Haar 
subspace $H$ satisfy {\rm (i')} and {\rm(ii)} of Corollary 3, for some integer $k \,(1 \le k \le n + 1)$. If 
$k$ elements $x_1, \ldots\hspace{-1pt}, x_k \in X$ are all distinct and the common value of {\rm(ii)}
is not zero, then there exists a positive number $\gamma$ such that }
\[ \max_{a \in A}\parallel f_a - h \parallel \ge \max_{a \in A}\parallel f_a - f^* \parallel
+ \, \gamma \parallel f^* - h \parallel~~for~all~h \in H.  \]
(Proof) First we show that $k = n + 1$. If $k \le n$, there is a function $f \in H$ such that
$f(x_i) = f_{a_i}(x_i) - f^*(x_i)$ for $i \,(1 \le i \le k)$, since $H$ is a Haar subspace. Inserting this $f$
into the equality (i') we have
\begin{eqnarray*} 
\sum_{i=1}^k  \lambda_i | f_{a_i}(x_i) -f^*(x_i) |^2 = | f_{a_i}(x_i) -f^*(x_i) |^2 = 0.  
\end{eqnarray*}
However, we assumed that $\delta = | f_{a_i}(x_i) -f^*(x_i)| = \max_{a \in A}\parallel f_a - f^*\parallel$ 
is not zero. Hence we must have $k = n + 1$. \\
\indent
Letting $h$ be an arbitrary element in $H$ such that $\parallel h\parallel = 1$, condition (i') can be
written as
\begin{eqnarray*} 
\sum_{i=1}^{n+1} \lambda_i \sigma_i h(x_i) = 0, ~~\sigma_i = (f_{a_i}(x_i) - f^*(x_i))/\delta~~(1 \le i \le n+1).
\end{eqnarray*}
Since $H$ is a Haar subspace and $\parallel h\parallel = 1$, it follows that 
$\max_{1 \le i \le n+1}\sigma_ih(x_i) > 0$.
If we set
\begin{eqnarray*} 
\gamma = \min_{\parallel h\parallel = 1} \,\max_{1 \le i \le n+1}\sigma_i h(x_i),
\end{eqnarray*}
then $\gamma$ is positive, since the set $\{h \in H : \; \parallel h\parallel =1 \}$ is compact. \\
\indent
Let $f \ne f^*$ be any element in $H$ and set $h = (f^* - f)/\parallel f^* - f\parallel$. Then
there exists at least one $i$ satisfying
\begin{eqnarray*} 
\sigma_i h(x_i) = \sigma_i (f^*(x_i) - f(x_i)) /\parallel f^* - f \parallel \ge \,\gamma,
\end{eqnarray*}
and the required inequality follows using this $i$ and $|\sigma_i| = 1$:
\begin{eqnarray*} 
\max_{a \in A}\parallel f_a -f \parallel \ge \sigma_i(f_{a_i}(x_i) - f(x_i))
= \sigma_i(f_{a_i}(x_i) - f^*(x_i)) + \sigma_i(f^*(x_i) - f(x_i)) \\
\ge \max_{a \in A}\parallel f_a - f^* \parallel + \,\gamma \parallel f^* - f \parallel.
\end{eqnarray*} 
\\
\\
{\large \bf 4. BSA on $L^p$-spaces} \\
\\
\indent
Let $(S, m)$ be a $\sigma$-finite positive measure space and $L^p(S, m)~(1 \le p < \infty)$ the set of all real-valued 
measurable functions $f$ such that $|f|^p$ are integrable over $S$. For such $p$
let $q$ be the real number determined by $p^{-1}+ q^{-1} = 1$ for $p>1$, and $q = \infty$ for $p = 1$. We use the following
notation (similarly for $\parallel g \parallel _q$ of $g \in L^q$), 
\begin{eqnarray*}
     \parallel f \parallel _p \, = \, \big(\int _S |f|^p dm \big)^{1/p}, 
\end{eqnarray*}
and $\parallel f \parallel _{\infty} \,=  {\rm ess ~sup}_{x \in S} |f(x)|$.  \\
\indent
We also assume that $A$ is a compact set of a Hausdorff topological space and to each $a \in A$ there corresponds a function $f_a$
which belongs to $L^p(S, m)$, and that the mapping $A \to L^p(S, m)$ so defined is continuous. \\
\indent
Let $H$ be an $n$-dimensional subspace of $L^p(S, m)$, where $n \ge 1$. The problem is to approximate simultaneously
the functions $\{f_a\}$ by elements of $H$. If $f^* \in H$ satisfies
\begin{eqnarray}
{\rm max}_{a \in A} \parallel f_a -f ^*  \parallel _p \; \le \;{\rm max}_{a \in A} \parallel f_a  - f \parallel _p
\end{eqnarray}
for all $f \in H$, we say that $f^*$ is a BSA to the functions $\{f_a\}$ from $H$.\\
\indent
In order to formulate this problem in relation to Lemma 1, we need the following well-known facts.
Property (a) is the Banach-Alaoglu theorem  (see [3, p.68]) and property (b) is the duality pairing (see [7, p.115]). 
\begin{itemize}
\item[{\rm (a)}] The dual space of $L^p(S, m)$ is $L^q(S, m)$ and  
$G = \{ g \in L^q(S, m) : \, \parallel g \parallel _q \, \le 1 \}$ is a compact set in the weak*-topology 
$\sigma(L^q(S, m), L^p(S, m))$.
\item[{\rm (b)}] For each $f \in L^p(S, m)$ we have $\parallel f \parallel _p \;=\, {\rm max}_{g \in G} \int _S (gf) dm$.
\end{itemize}
\indent
Therefore, (4) is equivalent to the following:
\begin{eqnarray}
\max_{(a, g) \in A\times G} \int _S g(f_a - f^*) dm \le \max_{(a, g) \in A\times G} \int _S g(f_a - f) dm
\end{eqnarray}
for all $f \in H$, and the problem is to find a function $f^*$ satisfying (5).\\ 
\indent 
It follows from (a) that $A \times G$ is a compact set in the product topology and it is easy to see that 
\[ J(f, a, g) = \int_S g(f_a - f) dm \]
is a jointly continuous function of the three variables $a \in A$, $g \in G$ and $f  \in H$, 
using H\"{o}lder's inequality and the definition of the weak*-topology. Moreover, it is a convex function 
with respect to $f$. Hence we can again invoke Lemma 1 for the characterization of best approximations.
\\
\\
{\bf Theorem 5.}  {\it An element $f^* \in H$ is a BSA to $\{f_a\}$ 
from $H$ if and only if, for some positive integer $k \,(1 \le k \le n + 1)$,
there exist $a_1, \ldots\hspace{-1pt}, a_k \in A$, $g_1, \ldots\hspace{-1pt}, g_k \in G$ and 
positive numbers $\lambda_1, \ldots\hspace{-1pt}, \lambda_k$ with sum one, satisfying the following two conditions}: \\
~~(i)~~\lefteqn {\int_S \big(\sum_{i=1}^k \lambda_i g_i \big)h \,dm = 0~~for~all~h\in H;}\\
~(ii)~~\lefteqn {\int_S g_i(f_{a_i}- f^*) dm = \, \parallel f _{a_i} - f^* \parallel _p \;
 = \max_{a \in A} \parallel f _a - f^* \parallel _p \,  ~for~all~i \,(1 \le i \le k).} \\
 \\
(Proof)  Let $f^*$ be a BSA. Define the set 
$U = \{ f \in H : {\parallel f - f^*\parallel}_p \le 1 \}$. Then 
$U$ is compact and convex, for $H$ is a finite-dimensional subspace of $L^p(S, m)$. 
The convexity follows from Minkowski's inequality (see [7, p.33]). It is obvious that $f^*$ also
minimizes $\max_{(a, g) \in A \times G}J(f, a, g)$ over the set $U$. It follows from Lemma 1,
for some $k\, (1 \le k \le n+1)$, 
that there exist $a_1, \ldots\hspace{-1pt}, a_k \in A$, $g_1, \ldots\hspace{-1pt}, g_k \in G$, and
numbers $\lambda_1, \ldots\hspace{-1pt}, \lambda_k > 0$ with $\sum_{i=1}^k \lambda_i=1$ such that
the following two inequalities hold: \\
\begin{eqnarray}
\sum_{i=1}^k\int _S \lambda_i g_i(f_{a_i}- f^*) dm 
\le \sum_{i=1}^k\int _S \lambda_i g_i(f_{a_i}- f) dm 
\end{eqnarray} 
for all $f \in U$; and
\begin{eqnarray}
\sum_{i=1}^{n+1}\int _S \mu_i h_i(f_{b_i}- f^*) dm 
\le \sum_{i=1}^k \int _S \lambda_i g_i(f_{a_i}- f^*) dm 
\end{eqnarray}
for all $(b_1, h_1), \ldots\hspace{-1pt}, (b_{n+1}, h_{n+1}) \in A \times G$, and all
nonnegative numbers $\mu_1, \ldots\hspace{-1pt}, \mu_{n+1}$ with sum one. \\
\indent
Inequality (6) implies
\begin{eqnarray*}
\sum_{i=1}^k \int _S \lambda_i g_i (f - f^*) dm \le 0 {~~\rm for ~all}\: f \in U. 
\end{eqnarray*}
Let us put $f = f^* + th$, where $h \in H$ is arbitrary and $t>0$ is so small that this $f$ belongs to $U$. 
Then we get
\begin{eqnarray*}
\sum_{i=1}^k \int _S \lambda_i g_ih \,dm = \int_S \big(\sum_{i=1}^k \lambda_i g_i \big)h \,dm
\le 0 {~~\rm for ~all}\: h \in H, 
\end{eqnarray*}
which means condition (i), since the left-hand side of the inequality must be zero. \\
\indent
By setting $b_i = a_i$ and $\mu_i = \lambda_i$ for all $i$~($1 \le i \le k$) in (7) and remarking
property (b), we see that (7) implies
\begin{eqnarray*}
\sum_{i=1}^k\lambda_i {\parallel f_{a_i}- f^*\parallel}_p 
\le \sum_{i=1}^k \int _S \lambda_i g_i(f_{a_i}- f^*) dm.
\end{eqnarray*}
Since the reverse inequality always holds, we conclude
\begin{eqnarray*}
\sum_{i=1}^k \lambda_i {\parallel f_{a_i}- f^*\parallel}_p  = 
\sum_{i=1}^k \int _S \lambda_i g_i(f_{a_i}- f^*) dm \le \max_{a \in A}\parallel f _a - f^* \parallel _p .
\end{eqnarray*}
Next in (7) putting $\mu_1 = 1$ (so $\mu_i = 0$ for $i \ne 1$) and $b_1 = a$ for any $a \in A$, 
we have for any $a \in A$
\begin{eqnarray*}
\max_{g \in G} \int _S g(f_a - f^*) dm = \; \parallel f _a - f^* \parallel _p \:
\le \sum_{i=1}^k \int _S \lambda_i g_i(f_{a_i}- f^*) dm, 
\end{eqnarray*}
hence
\begin{eqnarray}
\max_{a \in A} \parallel f _a - f^* \parallel _p \:
\le \sum_{i=1}^k \int _S \lambda_i g_i(f_{a_i}- f^*) dm \le \max_{a \in A} \parallel f _a - f^* \parallel _p. 
\end{eqnarray}
Then we conclude from (8) and $\lambda_i > 0$ that for all $i\,(1 \le i \le k)$
\begin{eqnarray*}
\parallel f _{a_i} - f^* \parallel _p \; 
= \int _S g_i(f_{a_i}- f^*) dm = \max_{a \in A} \parallel f _a - f^* \parallel _p \;,
\end{eqnarray*}
which is condition (ii). \\
\indent
Conversely, suppose that $f^*$ satisfies conditions (i) and (ii). 
Let $f \in H$ be any element. We have by (i)
\begin{eqnarray}
\sum_{i=1}^{k}\int _S \lambda_i g_i(f_{a_i}- f) dm - 
\sum_{i=1}^{k}\int _S \lambda_i g_i(f_{a_i}- f^*) dm  
= \int _S \big(\sum_{i=1}^{k}\lambda_i g_i \big) (f^*- f) dm = 0.
\end{eqnarray}
On the other hand, using (ii),
\begin{eqnarray*}
\sum_{i=1}^{k}\int _S \lambda_i g_i(f_{a_i}- f) dm - 
\sum_{i=1}^{k}\int _S \lambda_i g_i(f_{a_i}- f^*) dm  
\le \max_{1 \le i \le k} \parallel f _{a_i} - f \parallel _p  \\ - 
\max_{a \in A} \parallel f _a - f^* \parallel _p \; \le \max_{a \in A} \parallel f _a - f \parallel _p
- \max_{a \in A} \parallel f _a - f^* \parallel _p,
\end{eqnarray*}
where the condition $\sum_{i = 1}^k \lambda_i=1~(\lambda_i>0)$ is used. 
Immediately we conclude that $f^*$ is a BSA in view of (9), thereby completing the proof. \\
\\
\\
{\bf References}
\begin{itemize}
\item[{[1]}]  E. W. Cheney, {\it Introduction to Approximation Theory}, McGraw-Hill, 1966.
\item[{[2]}]  R. B. Holmes, {\it A Course on Optimization and Best Approximation}, Lecture Notes in
Mathematics, Vol. 257, Springer-Verlag, 1972.
\item[{[3]}]  W. Rudin, {\it Functional Analysis}~(2nd ed.), McGraw-Hill, 1991.
\item[{[4]}]  S. Tanimoto, Uniform approximation and a generalized minimax theorem, 
{\it Journal of Approximation Theory} {\bf 45}, 1--10, 1985.
\item[{[5]}]  S. Tanimoto, A characterization of best simultaneous approximations, {\it Journal of Approximation Theory} 
{\bf 59}, 359--361, 1989.
\item[{[6]}]  S. Tanimoto, On best simultaneous approximation, 
{\it Mathematica Japonica} {\bf 48}, 275--279, 1998.
\item[{[7]}]  K. Yosida, {\it Functional Analysis}~(3rd ed.), Springer-Verlag, 1971.
\end{itemize}
\end{document}